\documentclass[final,   nonatbib]{elsarticle}
\makeatletter
\let\c@author\relax
\makeatother
\usepackage{amssymb, amscd,  mathrsfs}
\usepackage{latexsym}
\usepackage{amsmath}
\usepackage{amsthm}
\usepackage{amsfonts, float}
\usepackage{lineno, hyperref, latexsym}
\usepackage{adjustbox}
\usepackage{subfig}
\usepackage{indentfirst}
\usepackage{enumerate,longtable, tabu}
\usepackage{mathtools}
\usepackage[backend=bibtex, style=ext-numeric, sortcites=true, isbn=false, maxnames=5, mincrossrefs=1, parentracker=true,giveninits=true, sorting=nyt,articlein=false]{biblatex}

\DeclareFieldFormat[article]{title}{{#1}\addperiod}
\AtBeginBibliography{}
\PassOptionsToPackage{unicode}{hyperref}
\DeclarePairedDelimiter{\diagfences}{(}{)}
\newcommand{\diag}{\operatorname{diag}\diagfences}

\graphicspath{{./img/}{./figures}}
\usepackage{tikz}

\usetikzlibrary{calc}

\usepackage[figurename=Fig.]{caption}

\modulolinenumbers[5]

\newtheorem{thm}{Theorem}[section]
\newtheorem{lem}{Lemma}[section]
\newtheorem{definition}{Definition}[section]

\newtheorem{proof*}{proof}[section]
\newtheorem{cor}{Corollary}[section]

\newtheorem{pro}{Proposition}[section]

\usepackage{bbding}

\providecommand{\sep}{;\;}

\makeatletter
\def\ps@pprintTitle{%
	\let\@oddhead\@empty
	\let\@evenhead\@empty
	\def\@oddfoot{\reset@font\hfil\thepage\hfil}
	\let\@evenfoot\@oddfoot}
\makeatother

\begin{document}

 \begin{frontmatter}

	\title{On the relationship between shortlex order and $A_\alpha$-spectral radii of graphs with  starlike branch tree}

	\author[mymainaddress]{Haiying Shan \corref{mycorrespondingauthor}}
	%	\cortext[mycorrespondingauthor]{Corresponding author}
	\ead{shan\_haiying@tongji.edu.cn}
	\address[mymainaddress]{School of Mathematical Sciences,   Tongji University,   Shanghai,   P. R. China}
	\author[muhuoaddress]{Muhuo Liu}
	\ead{liumuhuo@163.com}
	\address[muhuoaddress]{\small Department of Mathematics, South China Agricultural University, Guangzhou 510642, China}
	
	\date{}

	\begin{abstract}
Let $\mathcal{P}(n)$ denote the set of all partitions of $n$, whose elements are nondecreasing sequences of positive integers whose sum is $n$. For ${\bf a}=( n_{1}, n_{2},\ldots, n_{d}) \in \mathcal{P}(n)$, let $G({\bf a},v)$ denote the graph obtained from connected graph $G$ appending $d$ paths with lengths $n_{1},n_{2},\ldots,n_{d}$ on vertex $v$ of $G$.
We show that the ordering of graphs in $\mathcal{G}_{n}(v)=\{ G({\bf a},v)  \mid {\bf a} \in \mathcal{P}(n) \}$ by $A_\alpha$-spectral radii coincides with the shortlex ordering of $\mathcal{P}(n)$. 

	\end{abstract}

\end{frontmatter}

	\section{Basic Definitions}

For weighted graphs $G=(V,E,\omega)$  of order $n$, where $\omega: V \times V \rightarrow \mathbb{R}_{\geq 0}$ is the edge weight function such that  $\omega(u,v)> 0$  if and only if $uv \in E(G)$. The matrix $A(G) = (a_{ij})$ of order $n$ with $a_{ij}=w(v_iv_j)$ if $v_iv_j\in E(G)$ and $0$ otherwise  is called weighted adjacency matrix  of  $G$. As a special cases, $A(G)$ is equal to the  adjacency matrix  of $G$ when $w(e)=1$ for any edge $e \in E(G)$. For symmetric nonnegative square matrix $A=A(G)$,  the weighted graph $G$ is  denoted by $G(A)$. By the definition, the 
weighted graph $G(A)$ maybe have loops but no multiedges.

The weighted characteristic polynomial of $G$ is defined by
$\phi(G,x)=\det (xI-A(G)).$
The spectral radius $\rho(G)$ of a weighted graph $G$  is the largest eigenvalue  of its weighted adjacency matrix  $A(G)$.
 
If $H$ is a subgraph of $G$ with $w_{H}(e)=w_{G}(e)$ for each $e\in E(H)$, then $H$ is called a \textit{weighted subgraph} of $G$.

If $U\subset V$, then $G[U]$ is the subgraph of $G$ spanned by $U$. For $W\subset V(G)$, by $G-W$ we mean the subgraph $G[V-W]$. Similarly, for $F\subset E(G)$,   We also denote by $G-F$ the subgraph of $G$ obtained by deleting the edges of $F$. In particular,   we write shortly $G-e$, whenever $F=\{e\}$.

Coalescence of graphs is a significant operation involving two graphs. 
The coalescence of graphs can be naturally extended to weighted (di)graphs.
Let $G_1$ and $G_2$ be two disjoint weighted (di)graphs with $v_1\in V(G_1)$ and $v_2\in V(G_2)$. The $coalescence$ of $G_1$ and $G_2$, denoted by $G_1 (v_1)\cdot  G_2(v_2)$ (or $G_1 \cdot G_{2}$ for short), is obtained from $G_1$ and  $G_2$ by identifying $v_1$ with $v_2$ to form  a new vertex $u$ and the edge weighted function $\omega$  of  $G_1 (v_1)\cdot  G_2(v_2)$ will be defined as 
$$
\omega(e)= \begin{cases}
\omega_{G_{i}}(e),  & \text{if } e \in  E(G_{i})\setminus \{ (v_{i},v_{i}) \} \text{for } i=1,2; \\
\omega_{G_{1}}(v_{1},v_{1})+\omega_{G_{2}}(v_{2},v_{2})  & \text{if } e = (u,u); \\
0, & \hbox{otherwise.}
\end{cases}
$$

Let $G$ be a connected graph with root $v$  and denote by $G(\mathbf{a},v)$ the graph obtained from $G$ by appending $d$ paths with lengths $n_{1},n_{2},\ldots,n_{d}$ on $v$,  where $\mathbf{a}=( n_{1}, n_{2},\ldots, n_{d})$ is a nondecreasing sequence of  positive integers.  Especially, if $G$ is trivial, then  $G(\mathbf{a},v)$ is a starlike tree. Hereafter, we always write $S(\mathbf{a})$ or $S(n_{1}, n_{2},\ldots, n_{d})$  as a starlike tree, where $d>2$. 
We note that a starlike tree, often defined as a tree with exactly one vertex of degree greater than 2, can also be thought as the coalescence of at least three paths rooted at an end, which becomes its root. For convenience, we include a path rooted at any vertex as being a rooted starlike tree, with $1$ or $2$ paths.
It is clear that $G(\mathbf{a},v)=G \cdot S(\mathbf{a})$. In this case,  $S(\mathbf{a})$ is called a starlike branch tree of $G(\mathbf{a},v)$ at $v$.
Spectral properties of starlike trees have been reported in \cite{lepovic2001some,Oliveira_2018,stevanovic2022ordering}.

Let $G$ be a graph with adjacency matrix $A(G)$ and $D(G)$ be the diagonal matrix of its vertex degrees. 
For any real number  $\alpha\in [0,1)$, Nikiforov~\cite{MR3648656} proposed the problem to study the spectral properties of the family of matrices  $A_{\alpha}(G)$ defined as the convex linear combination:
\[
A_{\alpha}(G)=\alpha D(G)+ (1-\alpha)A(G).
\]
The spectral radius of $A_{\alpha}(G)$ is denoted by $\rho_\alpha(G)$.
For more information on the $A_{\alpha}$-spectra, we refer the reader to \cite{MR3786248, guo2018alpha1,Lin_2018} and their references.

\section{ Main Results}
Let $\mathcal{P}(n)$ denote the set of all partitions of $n$, whose elements are nondecreasing sequences of positive integers with sum equal to  $n$.  Let $\mathcal{P}(n, k)$  be the subset of $\mathcal{P}(n)$ whose  elements are all partitions of $n$ with length $k$.
Now we discuss shortlex order of partitions in $\mathcal{P}(n)$.   Shortlex ordering for $\mathcal{P}(n)$ is defined as follows:

\begin{definition}
For two partitions $\mathbf{a}=\left(a_{1},a_2, \ldots, a_{k}\right)$ and $\mathbf{b}=\left(b_{1},b_2, \ldots, b_p\right)$ in $\mathcal{P}(n)$, $\mathbf{a}$ precedes $\mathbf{b}$   if either $k<p$ or, when $k=p$, $a_{i}<b_{i}$ holds for the smallest index $i$ at which the two partitions differ, denoted by $\mathbf{a} \prec_{\mathrm{lex}} \mathbf{b}$. Clearly, shortlex ordering is a linear ordering.
\end{definition}

In \cite{Oliveira_2018}, Oliveira et al. obtained the following result:
\begin{thm}[\cite{Oliveira_2018}]\label{thm:starlike}
Let $\mathbf{a},\mathbf{b}$ be two arbitrary  partitions in $\mathcal{P}(n)$. If $S(\mathbf{a}) \ncong S(\mathbf{b})$, then
$\rho(S(\mathbf{a})) < \rho(S(\mathbf{b})) \iff   \mathbf{a} \prec_{\mathrm{lex}} \mathbf{b}.$
\end{thm}

The following is the main result of this paper:
\begin{thm}\label{mainthm}
	Let $G$ be a connected graph with root $u$ and  $\mathbf{a},\mathbf{b}$ be two arbitrary  partitions in $\mathcal{P}(n)$ with 
	$\mathbf{a} \prec_{\mathrm{lex}} \mathbf{b}$.  For \(\alpha \in [0,1)\), we have 
	\begin{equation}\label{mainequ}
		\rho_\alpha(G(\mathbf{a},u)) \leq \rho_\alpha(G(\mathbf{b},u)),
	\end{equation}
where the equality holds if and only if   $G(\mathbf{a},u) \cong G(\mathbf{b},u)$, that is, $G=K_1$ and the length of $\mathbf{b}$ is $2$.
%	where the equality holds if and only if $G=K_1$ and the length of $\mathbf{b}$ is $2$, that is, $G(\mathbf{a},u) \cong G(\mathbf{b},u)$ is a  path graph.
\end{thm}
It is easy to see that Theorem \ref{thm:starlike} is a special case of our main result for \(\alpha=0\) and $G$ being a trivial graph.  The method for  proof of  Theorem \ref{mainthm} is entirely different from that of Theorem \ref{thm:starlike} in paper \cite{Oliveira_2018}.

Before completing the proof of our main result, we shall introduce more notations, concepts, and useful tools that we shall use in the remainder of the paper.

Let $\mathbf{A}=\left(a_{i j}\right)$ be any real symmetric $n \times n$ matrix. Write $X=\{1,2,\ldots, n\},$ let $\pi = \left\{X_{1},X_2,\ldots, X_{m}\right\}$ be a partition of $X$ with $n_{i}=\left|X_{i}\right| $ for $1\le i\le m$.
Let  $A_{i j}=A[X_{i}:X_{j}]$ be the submatrix of  $A$ whose rows indexed by elements of $X_{i}$ and columns indexed by elements of $X_{j}$. If for any $1\le i\le j\le m$, each submatrix $A_{ij}$ have constant row sums, then  $\pi$ is said to be an {\em equitable partition} of $A$. Let $\mathbf{S}=\left(s_{i j}\right)$ be the $n \times m$ characteristic matrix of the partition, that is, $s_{i j}=1$ if $i \in X_{j},$ and 0 otherwise.      Then the matrices $Q_{l}(\mathbf{A})=\Lambda^{-1} S^{T} A S$  and $Q_{s}(\mathbf{A})=\Lambda^{-1 / 2} \mathbf{S}^{T} \mathbf{A} \mathbf{S} \Lambda^{-1 / 2}$ are called the {\em left quotient} and   
the {\em symmetric quotient} of $\mathbf{A}$, where $\Lambda=S^TS=\operatorname{diag}\left(n_{1},n_2, \ldots, n_{m}\right)$. Since $Q_{l}(A)$ and $Q_{s}(A)$ are similar, they have the same spectrum.

The weighted (di)graphs associated with $Q_{l}(A)$ and $Q_{s}(A)$ are called quotient (or divisor) graph  and 
symmetrized quotient graph of $A$ with respect to the partition $\pi$, respectively, denoted by $G(A)/\pi$ and $G_{s}(A)/\pi$. Equitable partitions and quotient graph represent a powerful tool in spectral graph theory.  We refer the reader to \cite{godsil2001algebraic,cvetkovic2010introduction, YOU201921} and the references therein
for properties and applications  of equitable partition and quotient graph.

The following result will play a key role in the proof of Theorem~\ref{mainthm}.

\begin{lem}[\cite{YOU201921,Atik_2019}]\label{lem:quotient}
	The spectral radius of a nonnegative square matrix $A$ is the same as the spectral radius of a quotient matrix corresponding to an equitable partition.
\end{lem}

% The following lemma is from the idea of L. Lov\'{a}sz and J. Pelik\'{a}n in \cite{lovasz1973eigenvalues}.
% \begin{lem}\label{lem:lovas}
% 	Let $G$ and $H$ be two weighted graphs of the same order. If $\phi(H,x) >\phi(G,x)$ for every $x$ in the interval $[\rho(G),\infty)$, we have $$\rho(H) < \rho(G).$$
% \end{lem}

\begin{lem}[\cite{BELARDO20101513}]\label{lem:HF}
Let $G u v H$ be the weighted graph obtained from $G$ and $H$ by adding a bridge uv of weight $\omega(u v)$, where $u \in V_{G}$ and $v \in V_{H}$. Then
$$
\phi(G u v H)=\phi(G) \phi(H)-\omega(u v)^{2} \phi(G-u) \phi(H-v).
$$	
\end{lem}

\begin{lem}[\cite{BELARDO20101513}]\label{swkink}
Let $G \cdot H$ be the coalescence of two rooted weighted digraphs $G$ and $H$ whose roots are $u$ and $v$, respectively. Then
$$
\phi(G \cdot H)=\phi(G) \phi(H-v)+\phi(G-u) \phi(H)-x \phi(G-u) \phi(H-v).
$$
\end{lem}

\begin{lem}[\cite{Minc1988nonnegative}]\label{lemc2}
The maximal eigenvalue of an irreducible matrix is greater than the maximal eigenvalue of its principal submatrices.
\end{lem}

\begin{lem}[\cite{Shan_2021}] \label{lemz2}
	Let $A$ and $B$ be two nonnegative matrices of order $n$ with $A \geq B$ and $A \neq B$. Then the following holds:
	\begin{center}
		$\phi(B,x) \geq \phi(A,x)$ for $x \geq \rho(A)$,
	\end{center}
	especially, when $A$ is irreducible, the inequality is strict.
\end{lem}

Suppose that $u$ is an end vertex of path graph $P_{n+1}$. Let $B_n$ be the principal submatrix of $A_\alpha(P_{n+1})$ obtained by deleting the row and column corresponding to the vertex $u$. Let $f_n$ be the characteristic polynomial of $B_n$ and $\theta_n=\rho(B_n)$. 
Here, we use the convention that the determinant of an empty matrix is $1$.  So $f_{0}=1$.
By Lemma \ref{lemc2}, we have $\theta_{n+1} > \theta_{n}$ for $n \geq 1$.

From the definition of $f_n$ and Lemma \ref{lem:HF}, the following proposition can be derived straightforwardly.

\begin{pro}\label{prop:path}
	Let $a,b, n$ be three positive integers  with $a+b=n$. Then for $\alpha \in [0,1)$, we have
	\begin{enumerate}[(1).]
		\item $f_0=1, f_1=x-\alpha$ and $f_n=(x-2\alpha)f_{n-1}-(1-\alpha)^{2}f_{n-2}$ for $n \geq 2$.
		\item $\phi_\alpha(P_n)=f_af_b-(1-\alpha)^{2}f_{a-1}f_{b-1}$. 
	\end{enumerate}
\end{pro}

\begin{lem}\label{lem:inequ}
	Let $a,b,l$ be three positive integers with $a=b+l$.  Then for $\alpha \in [0,1)$, the following holds when $x\geq \theta_l$:
	$$
	f_{a-1} f_b> f_{a}f_{b-1}
	$$
\end{lem}

\begin{proof}
	From (2) of Proposition \ref{prop:path}, we have
	\begin{align*}
		\phi_\alpha(P_{a+b-1})=f_af_{b-1}-(1-\alpha)^{2}f_{a-1}f_{b-2}
		=f_{a-1}f_{b}-(1-\alpha)^{2}f_{a-2}f_{b-1}.
	\end{align*}
	So %\vspace*{-10mm}
	\begin{align*}
		f_{a}f_{b-1}-f_{a-1}f_{b}=&(1-\alpha)^2(f_{a-1}f_{b-2}-f_{a-2}f_{b-1})\\
		=&(1-\alpha)^4(f_{a-2}f_{b-3}-f_{a-3}f_{b-2})=\cdots\\
		= & (1-\alpha)^{2(b-1)}(f_{0}f_{l+1}-f_{1}f_{l}) \\
		= & -(1-\alpha)^{2(b-1)}(\alpha f_{l}+(1-\alpha)^{2}f_{l-1}).
	\end{align*}
	Since $\theta_l$ and $\theta_{l-1}$ are the maximum real roots of $f_{l}$ and $f_{l-1}$, respectively and $\theta_{l}>\theta_{l-1}$,   $f_{l-1}>f_{l} \geq 0$ when $x\geq \theta_l$. So $f_bf_{a-1}-f_{b-1}f_{a}>0$.
\end{proof}

By taking Lemma \ref{swkink} in mind, the following result can be proven by using the same idea as   Lemma 3.1 of \cite{wang2021spectral} and so we omit its proof here.

\begin{lem}\label{lem:specomp}
	Let $G, H$ be two nontrivial weighted connected graphs with $ u \in V(G)$ and $v_{1}, v_{2} \in V(H)$. 
Take $G_{i}=$ $G(u)\cdot  H(v_{i})$ for $i=1,2$. If  $\phi\left(H-v_{2}\right)>\phi\left(H- v_{1}\right)$ for  $x \geq \rho\left(H-v_{1}\right)$, then    $\rho\left(G_{1}\right)<\rho\left(G_{2}\right)$.
\end{lem}

Let $v$ be a vertex of graph $G$ of order $m$ and graph $H=G(\mathbf{a},v)$ with $\mathbf{a}=[a]*s+[b]$, where $a,b,s$ are three integers with   $a>0, s>0$ and $b \geq 0$. This implies that  $H$ is the graph obtained from $G$ by appending $s$ paths of length $a$ and one path of length $b$ to the vertex  $v$ of $G$. Let $\pi$ be a partition of $V(H)$ with cells $C_1,C_2,\ldots, C_{m+a+b}$
where for $1\le i\le a$,  $C_i$ consists of all vertices on the $s$ pendent paths at $v$ of length $a$ at distance $i$ to $v$  and other cells are singletons. It is easy to see that $\pi$ is an equitable partition of $A_{\alpha}(H)$.   

Let $D'(G)$ be the diagonal matrix obtained from $D(G)$ by replacing the diagonal entry corresponding to vertices $v$ with $d_G(v)+s-1$. Take $Q'(G)=\alpha D'(G)+(1-\alpha)A(G)$ and denote the  weighted graph associated with $Q'(G)$ by $\widetilde{ G }$.
Notice that $\widetilde{ G }$ is trivial if only if $\alpha=0$ and $G$ is trivial.
Let $P_s(a,b)$ be weighted path graph on the vertex set $\{0,1,\ldots,a+b\}$ 
whose edge weights are $1$ except for $\omega(a-1,a)=\sqrt{s}$. Let $A(P_s(a,b))$ be the weighted adjacency matrix of $P_s(a,b)$ and $D(P_s(a,b))$ be the degree diagonal matrix of $P_s(a,b)$.
Take $Q(H)=\alpha D(P_s(a,b))+(1-\alpha)A(P_s(a,b))$ and denote the  weighted graph associated with $\alpha D(P_s(a,b))+(1-\alpha)A(P_s(a,b))$ by $\widetilde{ P}_s(a,b)$. Let $G_s(v,a,b)$ be the weighted graph obtained from $\widetilde{ G }$ and  $\widetilde{ P}_s(a,b)$ by coalescing vertices $v$ and $a$.
Then $G_s(v,a,b)$ is the symmetrized quotient graph of $A_{\alpha}(G({\bf a},v))$ corresponding to partition $\pi$.
Since $\pi$ is an equitable partition of $A_{\alpha}(G({\bf a},v))$, from Lemma~\ref{lem:quotient}, we have
\begin{equation}\label{equ:starlik}
	\rho_\alpha(G({\bf a},v))=\rho(G_s(v,a,b)).
\end{equation}
\begin{figure}[t]
	\vspace*{-1cm}
	\centering
	\subfloat[$H=G({\bf a},v)$]{\includegraphics[page=1,width=0.4\textwidth]{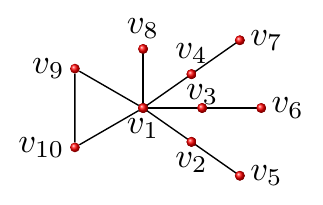}}\qquad
	\subfloat[$G_s(A_\alpha(H))/\pi$]{\vspace*{-1cm}\includegraphics[page=2,width=0.4\textwidth]{figure.pdf}}
	\caption{Graph $H$ and its symmetrized quotient graph of $A_\alpha(H)$}\label{fig:H}
	\vspace*{-4mm}
\end{figure}
To better understand the above notations, let us give an example.
Take ${\bf a}=[2,2,2,1]$, $G=C_3$ and $H=G({\bf a},v_{1})$(see (a) of Fig.~\ref{fig:H}).  Then $\pi =\{\{v_{2},v_{3},v_{4}\}, \{v_{5},v_{6},v_{7}\},\{v_{1}\},\{v_{8}\},\{v_{9}\},\{v_{10}\}\}$ 
is an equitable partition of $V(H)$ and $A_\alpha(H)$ whose rows and columns indexed by $V(H)$. $G_s(A_\alpha(H))/\pi$, the symmetrized quotient graph of $A_\alpha(H)$, is the associated graph of symmetrized quotient matrix of $A_\alpha(H)$ corresponding to $\pi$. The weight of loops at vertices $v_1,v_{2},v_{5},v_{8},v_{9},$ and $v_{10}$ are $6\alpha, 2\alpha,\alpha, \alpha, 2\alpha$ and $2\alpha$ respectively. And the weight of other edges in $G_s(A_\alpha(H))/\pi$ are $1-\alpha$ except for $\omega(v_{1},v_{2})=(1-\alpha)\sqrt{3}$.

For this example, $a=2, b=1$, $s=3$, $D(P_s(a,b))=\diag {1,2,2,1}$ and \\
$A(P_s(a,b))= \begin{pmatrix}
0 & 1 & 0 & 0 \\
1 & 0 & \sqrt{3} & 0\\
0 & \sqrt{3} & 0 & 1\\
0 & 0 & 1  & 0
\end{pmatrix}$, \quad
$Q'(G)= \begin{pmatrix}
	4 \alpha & 1-\alpha & 1-\alpha\\
	1- \alpha & 2 \alpha & 1-\alpha\\
	1 - \alpha & 1-\alpha & 2 \alpha
\end{pmatrix}.
$

For partition $\pi $, $Q_s(A_\alpha(H))$, the symmetric quotient matrix of $A_\alpha(H)$,
 is permutation-similar with the follow matrix
$$\begin{pmatrix}
	\alpha & 1-\alpha & 0 & 0 & 0 & 0\\
	1-\alpha & 2 \alpha & \sqrt{3}(1-\alpha) & 0 & 0 & 0\\
	0 & \sqrt{3}(1-\alpha) & 6\alpha & 1-\alpha & 1-\alpha & 1-\alpha\\
	0 & 0 & 1-\alpha  &  \alpha & 0 & 0\\
	0 & 0 &  1-\alpha & 0 & 2 \alpha & 1-\alpha\\
	0 & 0 & 1-\alpha  & 0 & 1-\alpha & 2 \alpha\\
	\end{pmatrix}.
	$$

For the sake of brevity, we define two operations of lists.
The concatenation of two lists ${\bf a}$ and ${\bf b}$ is denoted as ${\bf a}+{\bf b}$. Let ${\bf a}*n$ denote the concatenation of $n$ copies of list ${\bf a}$. For example, if ${\bf a}=[1,2,4]$ and ${\bf b}=[3,2]$, then ${\bf a} + {\bf b}=[1,2,4,3,2]$ and ${\bf a}*3=[1,2,4,1,2,4,1,2,4]$.

\begin{thm}\label{thm:quostarlik}
	Let $a,b,c,d$ and $s$ be five integers such that $a+b=c+d,$ $a > \max\{c,d\}\ge \min\{c,d\} > b \geq 0$ and $s>0$.  Take ${\bf e}=[a]*s+[b]$ and ${\bf f}=[c]*s+[d]$. Suppose that $u$ is a vertex of the  connected graph $G$. Then for $\alpha \in [0,1)$, the following holds
$$
\rho_{\alpha}(G({\bf e},u))\leq \rho_{\alpha}(G({\bf f},u)),
$$
where the equality holds if and only if $\alpha=0$ and $G$ is trivial graph.
\end{thm}
\begin{proof}
	From formula \eqref{equ:starlik}, we have 
	\begin{equation}\label{liu1e}\rho_\alpha(G({\bf e},u))=\rho(G_s(u,a,b)) \quad \text{and} \quad \rho_\alpha(G({\bf f},u))=\rho(G_s(u,c,d)).\end{equation}
	When $\alpha=0$ and $G$ is trivial graph, weighted graphs $G_s(u,a,b)$ and $G_s(u,c,d)$ are both isomorphic to $P_s(a,b)$. So $\rho_{\alpha}(G({\bf e},u))= \rho_{\alpha}(G({\bf f},u))$.
	
	When $\alpha>0$ or  $G$ is nontrivial, it is enough to show $$\rho(G_s(u,a,b))< \rho(G_s(u,a-1,b+1)).$$
	Let $H=\widetilde{P}_s(a,b)$ and rewrite vertices $a, a-1$ as $v_1, v_{2}$ respectively. Then $G_s(u,a,b)=\widetilde{G}(u)\cdot H(v_1), G_s(u,a-1,b+1)=\widetilde{G}(v) \cdot H(v_{2})$,

	$\phi(H-v_{1})= f_a f_b,$ and  $\phi(H-v_{2})= f_{a-1} f_{b+1}$. When    $x \geq \rho(G_s(u,a,b))> \theta_a$,  since $f_a f_b<f_{a-1} f_{b+1}$ by  Lemma \ref{lem:inequ}, we have $\phi(H-v_{1})<\phi(H-v_{2})$. From Lemma \ref{lem:specomp}, it follows that $\rho(G_s(u,a,b))<\rho(G_s(u,a-1,b+1))$. Now, from \eqref{liu1e}, we have  $\rho_\alpha(G({\bf e},u)) < \rho_\alpha(G({\bf f},u))$.\end{proof}

Let $G$ be a connected graph and $p,q$ be two nonnegative integers with $p \geq q+2$. Denote by $G_u(p,q)$  the graph obtained from $G$ by attaching two paths with lengths $p$ and $q$ to vertex $u$ of $G$. The results on the comparison of spectral radii of $G_u(p,q)$ and $G_u(p-1,q+1)$ can be found in \cite{ MR3786248, guo2018alpha1,qiao1979largest}.
Take $s=1, a=p,b=q$ in Theorem \ref{thm:quostarlik}, we will obtain the following result:
\begin{cor}[\cite{guo2018alpha1}]
	Let $G$ be a connected graph with $|E(G)| \geq 1$ and $u \in V(G)$. For integers $p \geq q \geq 1$ and $0 \leq \alpha<1, \rho_{\alpha}\left(G_{u}(p, q)\right)>\rho_{\alpha}\left(G_{u}(p+1, q-1)\right)$.
\end{cor}

\section{The proof of main result}

Let $\mathbf{a}, \mathbf{b}$  be two consecutive elements in  $\langle\mathcal{P}(n),\preceq_{\mathrm{lex}} \rangle$. Since $\langle\mathcal{P}(n),\preceq_{\mathrm{lex}} \rangle$ is linearly  order set, in order to prove Theorem \ref{mainthm} it suffices  to show $\rho_\alpha(G(\mathbf{a},u)) < \rho_\alpha(G(\mathbf{b},u))$ holds when $G(\mathbf{a},u) \ncong G(\mathbf{b},u)$.

Without loss of generality, suppose that $\mathbf{a}=\left(a_{1},a_2,\ldots, a_{k}\right) $ and $\mathbf{b}=(b_{1},b_2,$\par\noindent$\ldots, b_{l})$. Since $\mathbf{a} \prec_{\mathrm{lex}} \mathbf{b}$, $l\ge k$.  Actually, we only need to consider the  cases of $l=k$ or $l=k+1$.  For $k\le l\le k+1$, according to the definition of shortlex order, it is known that one of the following cases must occur  (see \cite{Oliveira_2018} or \cite{stevanovic2022ordering} for details):

\noindent {\bf Case I:} when $l=k$, there exists some $1\leq i \leq k-1$ such that $a_j=b_j$ for any $j <i$ and $b_{i}=b_{i+1}=\cdots =b_{k-1}=a_i+1$ and $\displaystyle b_k=a_k+\sum_{j=i}^{k-1}(a_j-a_i-1)$.
 
\noindent {\bf Case II:} when $l=k+1$,   $\mathbf{a}$ is  the maximum element of $\mathcal{P}(n,k)$ and $\mathbf{b}$ is the minimum element of $\mathcal{P}(n,k+1)$ in the linear order $\preceq_{\mathrm{lex}}$, as $\mathbf{a}$ and $\mathbf{b}$  are  two consecutive elements in  $\langle\mathcal{P}(n),\preceq_{\mathrm{lex}} \rangle$.  Suppose $n=qk+r$ with $0 \leq r < k$. Then\\
$
a_{1}=a_2=\cdots=a_{k-r}= \lfloor \frac{n}{k} \rfloor$ and $a_{k-r+1}=a_{k-r+2}=\cdots = a_{k}= \lceil \frac{n}{k} \rceil$, 
while
$$
b_{1}=b_2=\cdots=b_{k}=1, \quad \text{and} \,\,\,\, b_{k+1}=n-k \geq 1 .
$$
Since $n=qk+r \geq k \geq 1$ and $0 \leq r <k$, we have
\begin{equation}\label{Liu2e}n-k+1- \left\lceil \frac{n}{k} \right\rceil = (q-1)(k-1)+r -\left\lceil \frac{r}{k} \right\rceil \geq 0.\end{equation}

Since $G(\mathbf{a},u)\ncong G(\mathbf{b},u)$ by our hypothesis,  we may suppose that either   $G(\mathbf{a},u)$ or  $G(\mathbf{b},u)$ is not a path graph.

For Case II, take  $\mathbf{d}=[n-k+1]*k$. Then $G(\mathbf{a},u)$ is a proper subgraph of $G(\mathbf{d},u)$ by \eqref{Liu2e}. We have $\rho_\alpha(G(\mathbf{a},u)) < \rho_\alpha(G(\mathbf{d},u))$.  
So to prove Eq.\eqref{mainequ} for Case II, it is sufficient to show that  $\rho_\alpha(G(\mathbf{d},u))<\rho_\alpha(G(\mathbf{b},u))$ for $\mathbf{b} =[1]*k+[n-k]$.  It follows from Theorem \ref{thm:quostarlik} by taking $a=n-k+1, b=0, c=1, d=n-k$ and $s=k$.  Hence, $\rho_\alpha(G(\mathbf{a},u))<\rho_\alpha(G(\mathbf{b},u))$. 

For Case I, take $\widetilde{\mathbf{c}}=[a_1,a_2,\ldots,a_{i-1}],$ $\widetilde{\mathbf{a}}=[a_i,a_{i+1},  \ldots, a_k],$ and  $\widetilde{\mathbf{b}}=[a_i+1]*(k-i)+[b_k]$. Then, $\mathbf{a} =\widetilde{\mathbf{c}}+\widetilde{\mathbf{a}}$ and $\mathbf{b} =\widetilde{\mathbf{c}}+\widetilde{\mathbf{b}}$. We first show the following claim: \par\medskip\noindent {\bf Claim A.} $b_k+1 \geq a_k$ must hold. \par\medskip 
\noindent {\bf Proof of Claim A:} When $k-i=1$, from $a_i+a_{i+1}=b_i+b_{i+1}=a_i+1+b_k$, $b_k+1=a_k$,   Claim A holds.  For $k-i>1$, 
let us assume $b_k+1<a_k$ by contradiction. Then, $a_i+3\le a_k$ follows from  $b_k \geq b_{k-1}= a_i+1$. Note that $\sum_{j=i}^{k}a_{j}=\sum_{j=i}^{k}b_{j}=(k-i)(a_i+1)+b_k$ and $b_k-a_k\le -2$. Thus,  $\sum_{j=i+1}^{k-1}a_{j} < (a_i+1)(k-i-1)$. Combining this with  $\mathbf{a}$ being nondecreasing, we have $a_{i+1}=a_i$. Since  $a_{i+1}+3=a_i+3\le a_k$, we may let $
j_1$ be the maximum index  and $ j_2$ be the minimum index such that $a_{j_2}-a_{j_1}\ge 2$ and $i+1\le j_1<j_2\le k$.   By the choice of $j_1$ and $j_2$,   one of the following two situations must happen:
  \begin{enumerate}[(1).]
\item $j_2=j_1+1$, that is,    $a_{j_1+1}-a_{j_1}\geq 2$;
\item   $j_2\ge j_1+2$ and so   $a_{j_1}<a_{j_1+1}=\cdots=a_{j_2-1}<a_{j_2}$.
  \end{enumerate}
 For both situations, take 
 $\mathbf{p}=[p_1,p_2,\ldots,p_k] \in \mathcal{P}(n,k)$ such that $p_{j_1}=a_{j_1}+1, p_{j_2}=a_{j_2}-1$ and $p_{j}=a_{j}$ for $j \notin \{j_1, j_2\}$. Since $i+1\le j_1<j_2\le k$,  we have  
 $\mathbf{a} \prec_{\mathrm{lex}}  \mathbf{p}  \prec_{\mathrm{lex}}  \mathbf{b} $ for each situation, which contradicts the fact that $\mathbf{a} , \mathbf{b} $  are two consecutive elements in  $\langle\mathcal{P}(n),\preceq_{\mathrm{lex}} \rangle$. This confirms   Claim A. \qed
 
Let $G'=G(\widetilde{\mathbf{c}},u)$. Then $G(\mathbf{a},u)=G'(\widetilde{\mathbf{a}},u)$ and $G(\mathbf{b},u)=G'(\widetilde{\mathbf{b}},u)$.  
Let $\mathbf{c}=[a_i]+[b_k+1]*(k-i)$. Since $b_k+1 \geq a_k$ by Claim A, $G(\mathbf{a},u)$ is a subgraph of $G'(\mathbf{c},u)$.  We have $\rho_\alpha(G(\mathbf{a},u)) \leq \rho_\alpha(G'(\mathbf{c},u))$. 
Hence, 
  to prove Eq.\eqref{mainequ} for Case I, it is sufficient to show that  $\rho_\alpha(G'(\mathbf{c},u))<\rho_\alpha(G'(\widetilde{\mathbf{b}},u))$.  This  follows from Theorem \ref{thm:quostarlik} by taking $a=b_k+1, b=a_i, c=a_i+1, d=b_k$, $s=k-i$, and  $G=G'$, as $b_k\ge b_{k-1}=a_i+1>a_i$. So,   $\rho_\alpha(G(\mathbf{a},u))<\rho_\alpha(G(\mathbf{b},u))$.

This completes the proof of   Theorem \ref{mainthm}.

	\addcontentsline{toc}{section}{References}
	%\bibliography{ShanLiu}
	% \bibliography{110}
	\printbibliography
	
\end{document}